\documentclass[submission]{FPSAC2017}
\newtheorem{Thm}{Theorem}

\newtheorem{Prop}{Proposition}

\usepackage{graphicx,epsfig}
\title[Dimers, crystals and quantum Kostka numbers]{Dimers, crystals and quantum Kostka numbers}

\author{Christian Korff}

\address{School of Mathematics and Statistics, University of Glasgow\\
15 University Gardens, Glasgow G12 8QW, Scotland, UK}

\received{13 November 2016}

\abstract{We relate the counting of honeycomb dimer configurations on the cylinder to the counting of certain vertices in Kirillov-Reshetikhin crystal graphs. We show that these dimer configurations yield the quantum Kostka numbers of the small quantum cohomology ring of the Grassmannian, i.e. the expansion coefficients when multiplying a Schubert class repeatedly with different Chern classes. This allows one to derive sum rules for Gromov-Witten invariants.}


\begin{document}

\maketitle
%% note that you DO NOT have to put your abstract here -- it is generated by \maketitle and the \abstract and \resume commands above

%\section{Introduction}
%\label{sec:intro} 
\section{Dimer configurations on the cylinder}
This is an extended abstract of results for dimer configurations on the cylinder. Details and proofs will be provided in a separate publication.
The problem of counting dimer configurations on various lattices goes back to the works of Kasteleyn \cite{kasteleyn1967graph} and Fisher and Temperley \cite{temperley1961dimer}; see e.g. \cite{kenyon2009lectures} for an introduction and an overview. 

\begin{figure}[tbp]
\centering
\includegraphics[scale=0.5]{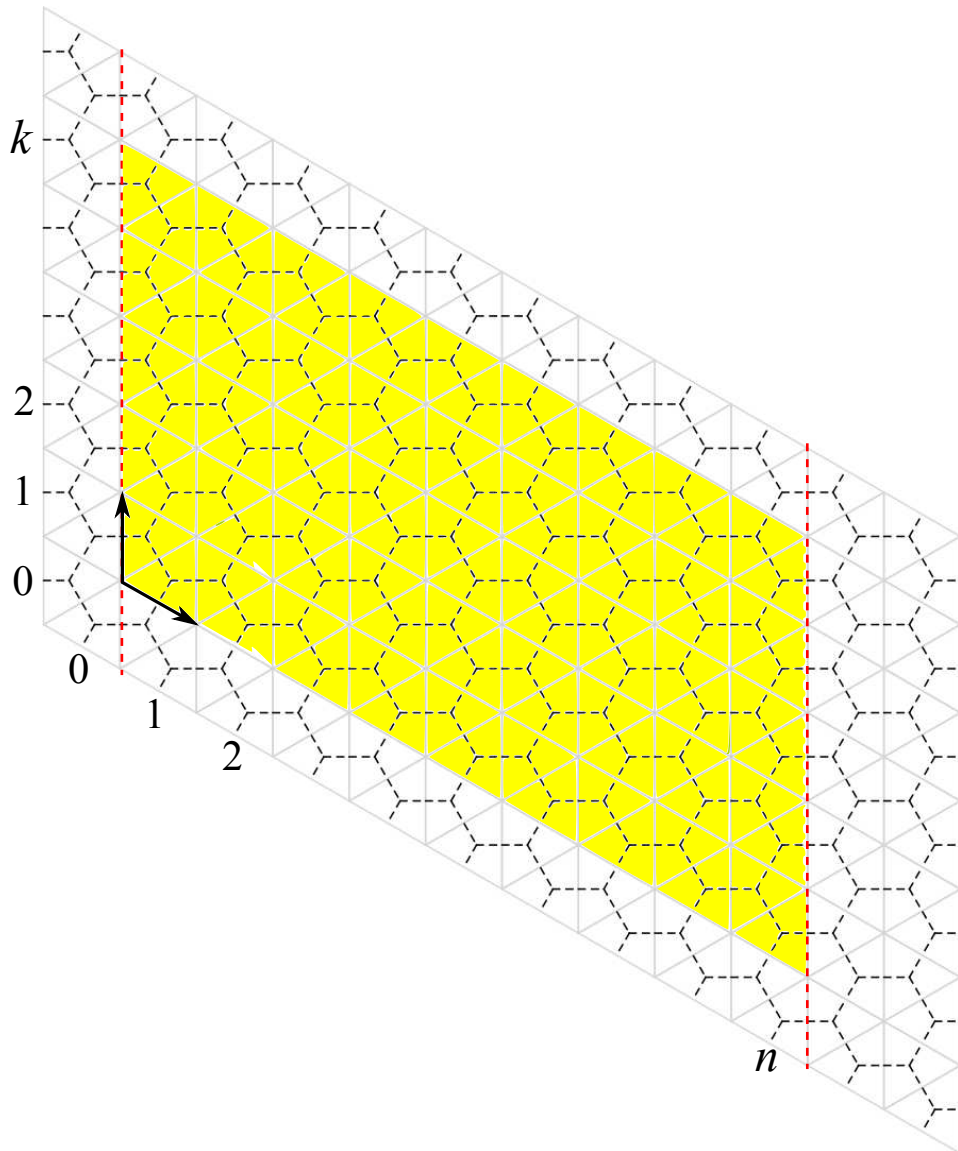}
\caption{Example of a honeycomb lattice on the cylinder when depicted in the plane. The yellow strip is glued together by identifying the lattice edges which are intersected by the red dotted lines, the boundary. The resulting cylinder has circumference $n=9$ and height $k=5$.}
\label{fig:honeycomb}
\end{figure}

\begin{figure}[tbp]
\centering
\includegraphics[scale=0.5]{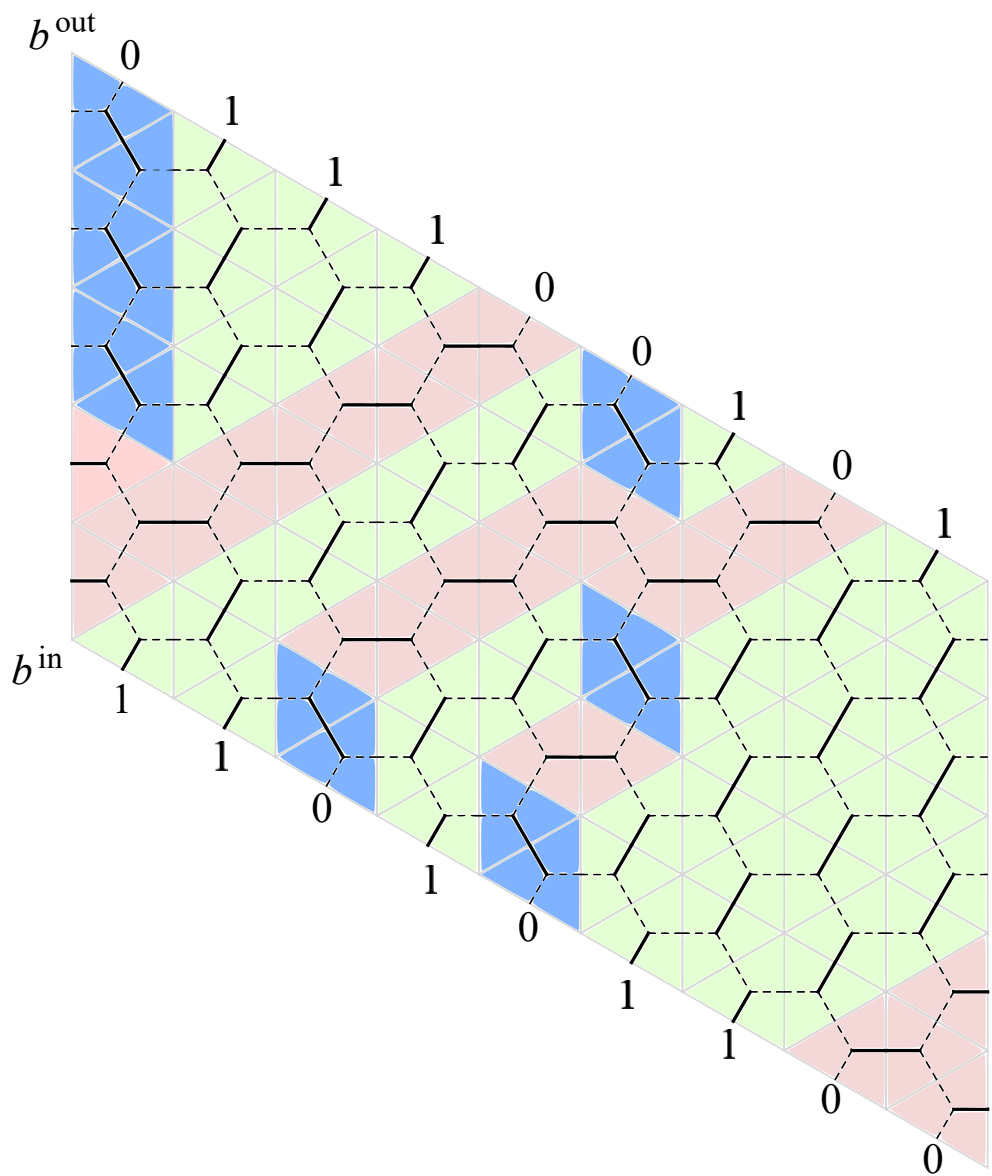}
\caption{Example of a dimer configuration on a cylinder with circumference $n=9$ and height $k=5$. The colouring of the lozenges
depicts the bijection to periodic tilings and plane partitions. The boundary conditions on the bottom and top of the cylinder are fixed by two binary strings $b^{in}$ and $b^{out}$ corresponding to the partitions $\mu=(4,4,3,2,2)$ and $\nu=(3,3,3,1,0)$. The number of horizontal dimers in each row (bottom to top) is given by $\lambda=(2,4,2,3,2)$ and there are $d=2$ dimers crossing the boundary in row 1 and 2.}
\label{fig:dimerconfig}
\end{figure}

Fix two integers $n\geq 3$ and $0\leq k\leq n$. Consider a hexagonal or honeycomb lattice on a cylinder of circumference $n$ and height $k$; see Figure \ref{fig:honeycomb} for an example. A {\em perfect matching} of the lattice is a selection of edges, called {\em dimers}, such that each vertex of the lattice is occupied by one and only one dimer; see Figure \ref{fig:dimerconfig} for an example. For simplicity, we shall visualise the cylindrical lattice as a strip of height $k$ in the Euclidean plane and consider only matchings with period $n$. We number the edges on the top and bottom of the strip consecutively with integers and call the vertical line separating the edges labeled with '0' and '1' the {\em boundary}; see Figure \ref{fig:honeycomb}. In what follows we are interested in the following refined counting problem of dimer configurations subject to a number of constraints or boundary conditions.

Firstly, we fix the dimer configurations on the top and bottom edges of the cylinder as shown in Figure \ref{fig:dimerconfig}, in terms of two binary strings $b^{in}=b_n^{in}\cdots b_2^{in}b_1^{in}$ and $b^{out}=b^{out}_n\cdots b_2^{in}b^{out}_1$, where a 1-letter signals a selected edge or dimer and a 0-letter an unoccupied edge. We require that $b^{in}$ and $b^{out}$ each contain $k$ one-letters, $\sum_{i=1}^n b^{in}_i=\sum_{i=1}^n b^{out}_i=k\,.$ %
Let $n\geq\ell^{in}_k>\cdots>\ell^{in}_1\geq 1$ and $n\geq\ell^{out}_k>\cdots>\ell^{out}_1\geq n$ be the positions of the 1-letters in $b^{in},b^{out}$ from right to left. Define two partitions $\mu,\nu$ by defining their parts through the relations $\mu_{k+1-i}+i=\ell^{in}_i$ and $\nu_{k+1-i}+i=\ell^{out}_i$, respectively. The Young diagrams of these partitions fit into a bounding box of height $k$ and width $k'=n-k$. The implicitly defined map between such partitions $\mu$ and binary strings of length $n$ with $k$ one-letters is a well-known bijection $\mu\mapsto b^\mu$ and, therefore, we shall identify both sets denoting them by the same symbol $\Pi_{k,n}$.

Secondly, we also fix the number of horizontal dimers in row $i$ to be $\lambda_i$ with $i=1,\ldots,k$ and set $|\lambda|=\sum_{i=1}^k\lambda_i$ to be the total number of horizontal dimers occurring in a configuration.

\begin{Thm} For given $\lambda,\mu,\nu\in\Pi_{k,n}$ denote by $\Gamma_{\lambda}(\mu,\nu)$ the set of perfect matchings or dimer configurations subject to the mentioned constraints with $b^{in}=b^\mu$ and $b^{out}=b^\nu$.
\begin{enumerate}
\item[(i)] If $\alpha$ is any permutation of the parts of $\lambda$, then  $|\Gamma_{\lambda}(\mu,\nu)|=|\Gamma_{\alpha}(\mu,\nu)|$.
\item[(ii)] The number of dimer configurations %
$|\Gamma_{\lambda}(\mu,\nu)|=0$ unless $|\lambda|+|\mu|-|\nu|$ is divisible by the circumference $n$ of the cylinder. 
\item[(iii)] If $|\Gamma_{\lambda}(\mu,\nu)|>0$ then in each perfect matching precisely
\begin{equation}\label{degree}
d=\frac{|\lambda|+|\mu|-|\nu|}{n}
\end{equation}
horizontal dimers cross the boundary.
\end{enumerate}
\end{Thm}

We can make a further statement about the possible minimum number of horizontal dimers if we only fix the boundary conditions $b^{in}$ and $b^{out}$ on the bottom and top of the cylinder but leave the number of horizontal dimers in each row arbitrary. For fixed $\mu ,\nu \in \Pi_{k,n}$ introduce the integers
\begin{equation}  \label{nmu}
n_i(\mu)=\sum_{j=n+1-i}^n b^\mu_j\,
\end{equation}
which are the partial sums of a binary string $b^\mu$ . Set 
\begin{equation}
d_{\min }(\nu ,\mu ):=\max_{i\in [n]}\{n_{i}(\nu )-n_{i}(\mu )\}\,  \label{dmin}
\end{equation}
and denote by $\Gamma(\mu,\nu)=\cup_\alpha\Gamma_{\alpha}(\mu,\nu)$ with $\alpha$ now a composition of length $\leq k$ with $\alpha_i\leq k'=n-k$. Then we have the following:
\begin{Prop}\label{mindimer}
(i) The minimal number of horizontal dimers in any perfect matching $\gamma\in\Gamma(\mu,\nu)$ is given by $|\lambda^{\min}|=n d_{\min}+|\nu|-|\mu|$ and in that configuration precisely $d_{\min}$ dimers are crossing the boundary.  
(ii) If $\lambda\in\Pi_{k,n}$ in \eqref{degree} is such that $d>d_{\min }(\nu ,\mu )$, then $|\Gamma_\lambda(\mu,\nu)|>0$ if and only if 
\begin{equation}
\theta _{i}(\mu ,\nu ;d):=d+n_{i}(\nu)-n_{i}(\mu )>0,\quad \quad 1\leq
i\leq n\;.  \label{theta}
\end{equation}%
The last constraint translates into the requirement that in each column of the lattice there is at least one horizontal dimer.
\end{Prop}

\begin{figure}[tbp]
\centering
\includegraphics[scale=0.5]{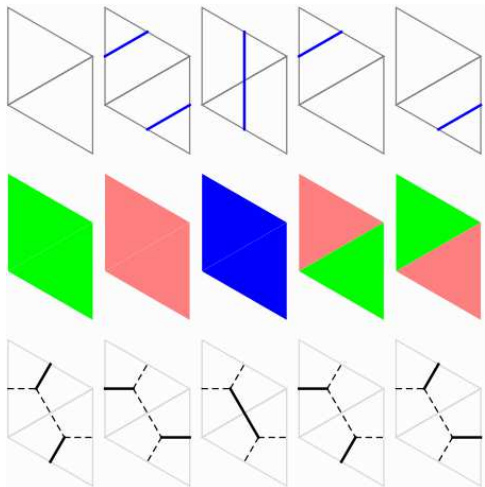}
\caption{Consider tilings of the plane using the depicted lozenges in each row, such that adjacent tilings match. The bottom row yields the dimer configurations from Figure \ref{fig:dimerconfig} if lozenges are placed such that dimers connect to dimers. Replacing each dimer lozenge with the corresponding one in the row above we obtain plane partitions and if we use the top row then we obtain domain walls or nonintersecting paths; see Figure \ref{fig:dimerconfig} for an example of how to map to a plane partition with periodic boundary conditions.}
\label{fig:rosetta}
\end{figure}

There are known bijections between the discussed dimer configurations and plane partitions or lozenge tilings as well as domain walls (non-intersecting paths) describing the spin configurations of the ground state of the triangular antiferromagnetic Ising model; see Figure \ref{fig:rosetta}. The stated results can therefore be reformulated in terms of any of these combinatorial tilings. It is the domain wall picture, see Figure \ref{fig:rosetta}, which was discussed in \cite{korff2014quantum} in connection with the small quantum cohomology of the Grassmannian. One then recognises the minimal number of horizontal dimers \eqref{dmin} as the minimal power of the deformation parameter $q$ occurring in a product of the small quantum cohomology ring Grassmannian; compare with \cite{yong2003degree} and \cite{fulton2004quantum}.

%\begin{figure}[tbp]
%\centering
%\includegraphics[scale=0.5]{PP2DW.eps}
%\caption{The mapping of the lozenge tiling from Figure \ref{fig:dimerconfig} onto non-intersecting paths corresponding to domain walls
%(shown in yellow) in the groundstate of the triangular antiferromagnetic Ising model.}
%\label{fig:DWconfig}
%\end{figure}

\section{Kirillov-Reshetikhin crystals}

Kashiwara's crystal bases \cite{kashiwara1991crystal} and their associated coloured, directed graphs, called {\em crystal graphs} or simply {\em crystals}, 
are an important combinatorial tool in representation theory; see e.g. \cite%
{hong2002introduction} for a textbook and references therein. A crystal graph consists of a set of vertices $B$, the basis elements, 
and certain maps $e_{i},f_{i}:B\rightarrow B\sqcup \{\emptyset \}$, $i\in [n]$ called 
\emph{Kashiwara operators} which define the directed, coloured edges of the graph: there exists 
an edge $b\rightarrow b'$ of colour $i $, if and only if, $f_{i}(b)=b'$ in which case we also must have $e_{i}(b')=b$. 
In particular, there are no multiple edges.

Given a crystal graph and a vertex $b\in B$ one can consider the maximal length
of a directed path along edges of a fixed colour $i$ which ends or starts at $b$, 
\begin{equation}
\varepsilon _{i}(b)=\max \{p\in \mathbb{Z}_{\geq 0}:e_{i}^{p}(b)\neq
\emptyset \},\quad \varphi _{i}(b)=\max \{p\in \mathbb{Z}_{\geq
0}:f_{i}^{p}(b)\neq \emptyset \}\;.  \label{eps&phi}
\end{equation}%
The functions  \eqref{eps&phi} allow one to
introduce the tensor product $B_{1}\otimes B_{2}$ of two crystal graphs $B_1,B_2$ as the crystal
obtained through the following action of the Kashiwara operators $e_i,f_i:%
B_{1}\times B_{1}\rightarrow B_{1}\times B_{2}\sqcup
\{\emptyset \}$ on the Cartesian product of the respective vertex sets, 
\begin{eqnarray}
e_{i}(b_{1}\otimes b_{2}) &=&\Bigl\{%
\begin{array}{ll}
e_{i}(b_{1})\otimes b_{2}, & \varepsilon _{i}(b_{1})>\varphi _{i}(b_{2}) \\ 
b_{1}\otimes e_{i}(b_{2}), & \text{else}%
\end{array}%
\Bigr.  \label{ctensore} \\
f_{i}(b_{1}\otimes b_{2}) &=&\Bigl\{%
\begin{array}{ll}
f_{i}(b_{1})\otimes b_{2}, & \varepsilon _{i}(b_{1})\geq \varphi _{i}(b_{2})
\\ 
b_{1}\otimes f_{i}(b_{2}), & \text{else}%
\end{array}%
\Bigr.  \label{ctensorf}
\end{eqnarray}%
together with the convention $b\otimes \emptyset =\emptyset $ and $\emptyset
\otimes b=\emptyset $. Note that
there exist different conventions for the definition of the tensor product
in the literature, our choice will be suited for the discussion at hand.

Here we are concerned with tensor products
\begin{equation}
B_{\lambda }:=B_{\lambda _{1}}\otimes \cdots \otimes B_{\lambda _{k}}\,,\qquad\lambda\in\Pi_{k,n},
\end{equation}%
of the crystal graphs $B_{r},\;0\leq r<n$ (usually denoted by $B^{r,1}$ in the
literature) of certain finite-dimensional level one modules of the affine
quantum algebra $U_{\upsilon }(\mathfrak{\widehat{sl}}_{n})$, so-called
Kirillov-Reshetikhin (KR) modules \cite{chari1995quantum}. The basis elements in $B_r$ are 
labelled by binary strings $b$ of length $n$ with $r$-one letters, hence as sets we have $B_r=\Pi_{r,n}$. The basis elements in 
$B_\lambda$ are then $k$-tuples $%
b^{(1)}\otimes \cdots \otimes b^{(k)}$ of binary strings $b^{(j)}$ which can be efficiently 
recorded in terms column tableaux, where column $j$ contains the positions 
$\ell^{(j)}_i$ of 1-letters of $b^{(j)}$. 

Not every module of $U_{\upsilon }(\mathfrak{%
\widehat{sl}}_{n})$ possesses a crystal basis, but KR modules are
distinguished by the fact that they do and that the associated crystal
graphs are \emph{perfect}. That is, tensor products of KR modules
again possess crystal bases and their associated crystal graphs are
connected; see \cite{hong2002introduction} and references therein.

%\begin{figure}[tbp]
%%\begin{equation*}
%\centering
%\includegraphics[scale=0.9]{KRcrystal_2_7.eps}
%%\end{equation*}%
%\caption{Example of a single column KR crystal for $N=7$ and $r=2$. The
%outer vertices of the star shaped crystal graph correspond to 01-words with
%adjacent 1-letters. For instance, the lowest vertex is the 01-word $0110000$%
%. }
%\label{fig:KRcrystal}
%\end{figure}

\subsection{Combinatorial R-matrix and dimers}

Fix $\lambda\in\Pi_{k,n}$. There exists a unique
bijection $R_{\lambda}: B_{n-k}\otimes B_{\lambda }\rightarrow
B_{\lambda}\otimes B_{n-k}$, called \emph{combinatorial $R$-matrix}, which
preserves the crystal graph structure and is a (set-theoretical) solution of the Yang-Baxter equation; see e.g. \cite{nakayashiki1997kostka}. 
In addition, the $\widehat{\mathfrak{sl}}_{n}$ Dynkin diagram automorphism $\Omega$ induces another trivial graph isomorphism $%
\Omega: B_{\lambda}\rightarrow B_{\lambda}$ by cyclic permutations of the letters in the binary strings, i.e. $\Omega(b^{(1)}\otimes \cdots \otimes b^{(k)})=\Omega(b^{(1)})\otimes \cdots \otimes \Omega(b^{(k)})$ with $\Omega(b^{(j)})=b_1^{(j)}b_k^{(j)}b_{k-1}^{(j)}\cdots b_2^{(j)}$. It commutes with the combinatorial $R$-matrix. We now define a particular subset of crystal vertices $b\in B_\lambda$.

\begin{figure}[tbp]
%\begin{equation*}
\centering
\includegraphics[scale=0.5]{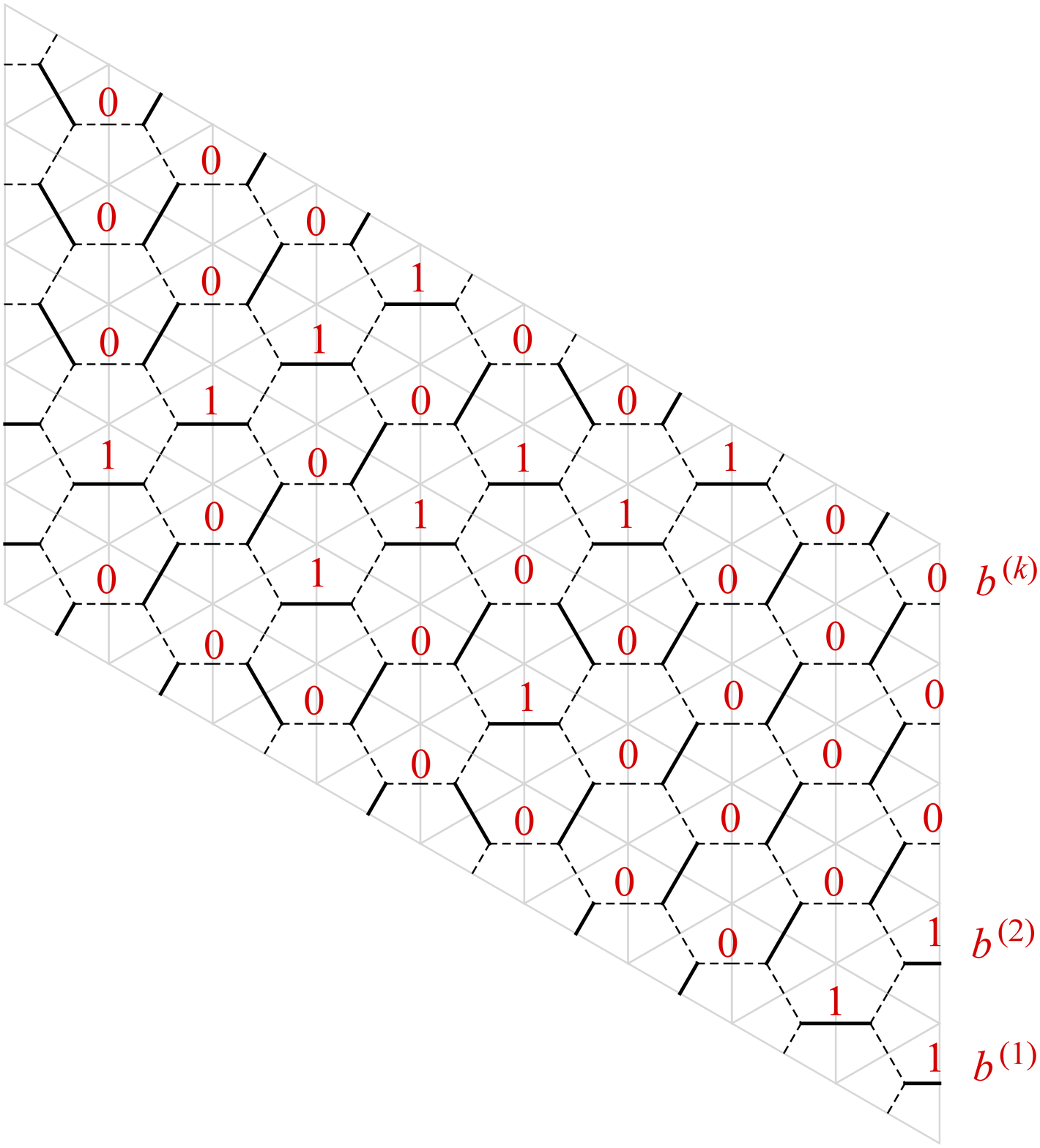}
%\end{equation*}%
\caption{Example of mapping the dimer configuration from Figure \ref{fig:dimerconfig} onto a vertex in the KR crystal %
$B_{2,4,2,3,2}=B_2\otimes B_4\otimes B_2\otimes B_3\otimes B_2$ by recording in each lattice row the occurrence of a horizontal dimer with a 1-letter and otherwise putting down a 0-letter. This results in a $k=5$-tuple $b^{(1)}\otimes\cdots\otimes b^{(k)}$ of binary strings which satisfy the conditions from Prop \ref{prop:crystalset}.}
\label{fig:dimer2KRcrystal}
\end{figure}

For fixed $\nu,\mu\in\Pi_{k,n}$ set 
\begin{equation*}
\Theta _{i}=\left\{ 
\begin{array}{ll}
1, & \theta _{i}(\mu ,\nu ,d)=d+n_{i}(\nu )-n_{i}(\mu )>0 \\ 
0, & \text{else}%
\end{array}%
\right. \;,
\end{equation*}
where $\theta$ is the same integer vector as in \eqref{theta}. Denote by $\mu',\nu'\in\Pi_{n-k,n}$ the 
conjugate partitions of $\mu,\nu$ obtained by transposing the respective Young diagrams.
\begin{Prop}\label{prop:crystalset}
Let $b\in B_{\lambda}$. The following
statements are equivalent.

\begin{itemize}

\item[(i)] $\varphi_i(b)=b^{\mu }_{n-i}\Theta _{i}$\quad and\quad $\varepsilon_i
(b)=b^{\nu }_{n-i}\Theta _{i+1}$ with $i\in\mathbb{Z}_n$.

\item[(ii)] $R_{\lambda}(b^{\mu' }\otimes \Omega%
(b))=b\otimes b^{\nu'}$.
\end{itemize}
 \end{Prop} 
Suppose $d>d_{\min }$. Then property (i) simplifies to%

\begin{equation}
\varphi_i(b)=b^{\mu }_{n-i}\quad \text{and}\quad \varepsilon_i(b)=b^{\nu
}_{n-i}\;.
\end{equation}
This characterisation of crystal vertices is an affine extension of the one considered by 
Berenstein and Zelevinsky in \cite{berenstein2001tensor} when describing Kostka numbers and Littlewood-Richardson coefficients for type $A$. Their results extend to all finite semi-simple Lie algebras.
\begin{Thm}
Denote by $B_{\lambda}(\mu,\nu)$ the set of crystal graph vertices $b\in B_\lambda$ satisfying the conditions of the previous proposition. There exists a bijective map between the elements in $\Gamma_\lambda(\mu,\nu)$ and $B_\lambda(\mu,\nu)$; see Figure \ref{fig:dimer2KRcrystal} for an example.
\end{Thm}
In other words, the $i$-signatures of the elements in $B_{\lambda}(\nu ,\mu )$ are fixed in terms of the start and end positions of the
dimer configurations and any crystal vertex with these $i$-signatures must be
the image of such a dimer configuration. 

\section{Quantum cohomology and toric Schur functions}
Quantum cohomology arose from works of Gepner \cite{gepner1991fusion},
Intriligator \cite{intriligator1991fusion}, Vafa \cite{vafa1991topological}, Witten \cite{witten1993verlinde} and since then has been studied extensively in the mathematics literature. The small quantum cohomology ring of the Grassmannian of $%
k $-planes in $\mathbb{C}^{n}$ has the following known presentation \cite%
{siebert1997quantum} 
\begin{equation}
qH^{\ast }(\text{Gr}_{k}(\mathbb{C}^n))\cong \mathbb{Z}[q][e_{1},\ldots
,e_{k}]/\langle h_{n-k+1},\ldots ,h_{n-1},h_{n}+q(-1)^{k}\rangle ,
\label{qH*}
\end{equation}%
where $h_{r}=\det (e_{1-i+j})_{1\leq i,j\leq r}$ are the Chern classes of the normal vector bundle and the $e_i$'s are the Chern classes of the canonical bundle. Denote by  $\sigma_{\lambda }=\det (h_{\lambda
_{i}-i+j})_{1\leq i,j\leq k}$ the Schubert classes and consider the 
coefficients in the following product expansion in (\ref{qH*}), 
\begin{equation}
\sigma_{\mu }\ast h_{\lambda _{1}}\ast \cdots \ast h_{\lambda _{r}}=\sum_{\nu \in\Pi_{k,n}}q^{d}\sigma_{\nu }K_{\nu /d/\mu ,\lambda },  \label{qKostkadef}
\end{equation}%
which are called \emph{quantum Kostka numbers} \cite{bertram1999quantum}. Here the power $d$,
called `degree', is fixed through the equality \eqref{degree}%
, otherwise the coefficient vanishes. 

\begin{Thm}[Sum rule for quantum Kostka numbers]
\label{sumruleK} 
(i) The number of possible dimer configurations
and crystal vertices fixed by the partitions $\lambda,\mu,\nu\in\Pi_{k,n}$ matches the quantum Kostka
number, 
\begin{equation}
K_{\nu/d/\mu,\lambda}=|\Gamma_{\lambda}(\mu,\nu)|=|B_{\lambda}(\mu,\nu)|\;.
\end{equation}
(ii) Summing over all compositions $\alpha=(\alpha_1,\ldots,\alpha_k)$ with $0\leq\alpha_i\leq n-k$, one obtains the 
total number $|\Gamma(\mu,\nu)|$ of dimer configurations on the cylinder subject only to the boundary conditions 
$b^{in}=b^\mu$ and $b^{out}=b^\nu$ on the bottom and top of the cylinder, 
\begin{equation}
|\Gamma(\mu,\nu)|=\sum_{\alpha}|\Gamma_{\alpha}(\mu,\nu)|= 
\sum_{\lambda\in\Pi_{k,n}}K_{\nu/d/\mu,\lambda}
\frac{\ell(\lambda)!}{\prod_{i\geq 1}m_i(\lambda)!}\binom{k}{\ell(\lambda)}\,,
\end{equation}%
where $\ell(\lambda)$ is the length of the partition $\lambda$ and $m_i(\lambda)$ the multiplicity of $i$ in $\lambda$.
\end{Thm}

\subsection{Toric Schur functions}
The predominant mathematical interest in the ring (\ref{qH*}) is the
computation of the 3-point genus 0 Gromov-Witten invariants $C_{\lambda \mu
}^{\nu ,d}$. The latter occur in the product expansion of two Schubert
classes %
%$s_{\lambda }$ in $qH^{\ast }(%
%\limfunc{Gr}\nolimits_{n,n+k})$, 
\begin{equation}
\sigma_{\mu }\ast \sigma_{\lambda }=\sum_{d\geq 0,\nu \in (n,k)}q^{d}C_{\lambda \mu
}^{\nu ,d}\sigma_{\nu },
%\quad \sigma_{\lambda }=\det (e_{\lambda _{i}^{\prime}-i+j})_{1\leq i,j\leq k},  
\label{GWinv}
\end{equation}%
and count rational curves of degree $d$ intersecting three Schubert
varieties in general position, which are parametrised by $\lambda ,\mu ,\nu $%
; for details we refer the reader to the literature, e.g. \cite{bertram1997quantum}%
, \cite{bertram1999quantum}, \cite{buch2003quantum}, \cite{buch2003gromov} and references therein.

A combinatorial interpretation of Gromov-Witten invariants was given in \cite{postnikov2005affine}: one generalises the 
notion of an ordinary skew Schur function $s_{\nu/\mu}$, where the expansion coefficients in the basis of 
Schur functions $s_\lambda$ are given by Littlewood-Richardson coefficients, $%
c_{\lambda \mu }^{\nu }=C_{\lambda \mu }^{\nu ,0}$, to so-called toric Schur 
functions,%
\begin{equation}  \label{toricschur}
s_{\nu /d/\mu }(x_{1},\ldots ,x_{k})=\sum_{\lambda \in\Pi_{k,n}}C_{\lambda \mu
}^{\nu ,d}s_{\lambda }(x_{1},\ldots ,x_{k})\;.
\end{equation}
Postnikov introduced these functions in terms of so-called toric tableaux \cite{postnikov2005affine}, which are special cases of the cylindric plane partitions considered by Gessel and Krattenthaler in \cite{Gessel1997cylindric}. Here we express them as weighted sums over KR crystals and dimer configurations.
\begin{Prop}
Toric Schur functions can be expressed as the following weighted sums,   
\begin{equation}
s_{\nu /d/\mu }(x_{1},\ldots ,x_{k})=\sum_{\alpha }| B_{\alpha}(\mu,\nu)|x^{\alpha }
=\sum_{\alpha }|\Gamma_{\alpha}(\mu,\nu)|x^{\alpha }\,,
\end{equation}%
where  $\alpha$ runs over all \emph{compositions} which have at
most $k$ parts $\alpha_i\leq n-k$.
\end{Prop}
If the degree \eqref{degree} vanishes, $d=0$,
one has $s_{\nu /0/\mu }=s_{\nu /\mu }$ and in this case one can use the
Robinson-Schensted-Knuth correspondence to arrive at the familiar crystal
theoretic interpretation of skew Schur functions.  This interpretation can be extended to $d=d_{min}$ using the cyclic $\mathbb{Z}_n$ symmetry of the cylinder manifest in Prop \ref{prop:crystalset} (ii), similar to the discussion in \cite{morse2012combinatorial}.

From the expansion \eqref{toricschur} one now arrives at the following:
\begin{Thm}[Sum rule for Gromov-Witten invariants]
One has the following alternative sum rule for the total number $|\Gamma(\mu,\nu)|$ of dimer configurations on the cylinder,
\begin{equation}
|\Gamma(\mu,\nu)|=\sum_{\alpha}|\Gamma_{\alpha}(\mu,\nu)|= \sum_{\lambda\in\Pi_{k,n}}C_{\lambda
\mu }^{\nu ,d}\prod_{s\in \lambda }\frac{n+c(s)}{h(s)},
\label{sumrule0}
\end{equation}%
where the product runs over all squares $s=\langle i,j\rangle $ in the Young
diagram of $\lambda$ and $c(s)=j-i$ denotes its content and $%
h(s)=\lambda_i+\lambda\rq{}_j-i-j+1$ its hook length.
\end{Thm}
A different connection between Gromov-Witten invariants and (combinatorially defined) crystals has been found by Morse and Schilling in \cite{morse2016crystal}. In {\em loc. cit.} the authors define a crystal structure on particular factorisations of affine permutations and identify the Gromov-Witten invariants of the full flag variety with the number of certain highest weight factorisations of affine permutations. They recover the Gromov-Witten invariants for the Grassmannian as special case of their more general construction \cite[Thm 5.16]{morse2016crystal}.

In contrast the results here connect the small quantum cohomology ring with the known crystal structure of KR modules of $U_{\upsilon }(\mathfrak{\widehat{sl}}_{n})$ and their combinatorial $R$-matrix. We hope to make the connection between both crystal structures in future work. We believe that such a connection would help with the construction of `quantum group structures', so-called Yang-Baxter algebras which provide maps $qH^{\ast }(\text{Gr}_{k}(\mathbb{C}^n))\rightarrow qH^{\ast }(\text{Gr}_{k\pm 1}(\mathbb{C}^n))$ \cite{korff2010wznw,korff2014quantum}, to general flag varieties; see also \cite{gorbounov2014equivariant} for an extension of the discussion to equivariant quantum cohomology and \cite{gorbounov2014quantum} for quantum K-theory.

\acknowledgements{
Part of this research was started back in 2011 while visiting the Hausdorff Institute for Mathematics (HIM), Bonn and while being on a University Research Fellowship of the Royal Society. The author is also indebted to Anne Schilling and Catharina Stroppel for discussions and sharing knowledge.
}

%% if you use biblatex then this generates the bibliography
%% if you use some other method then remove this and do it your own way
%\printbibliography

\end{document}